\documentclass[a4paper,12pt]{article}
\usepackage{amssymb}
\usepackage{amsthm}
\usepackage{mathrsfs}

\newtheorem{thm}{Theorem}[section]
\newtheorem{prop}[thm]{Proposition}
\newtheorem{cor}[thm]{Corollary}
\newtheorem{lm}[thm]{Lemma}
\newtheorem{prob}[thm]{Problem}
\newtheorem{rmk}[thm]{Remark}

\def\hline{\smash{\mathop{\hbox to 5mm{\hrulefill}}}}

\def\References #1\par{\vskip24pt plus 1pt
\setbox0=\hbox
{\ninepoint [#1]\enspace}
\centerline{\twelverm References}\nobreak\par\nobreak
\interlinepenalty=10000\parskip=0pt plus 1pt
\parindent=\wd0
\ninepoint
\nobreak\vskip12pt\nobreak}
\def\refer #1 {\litem{[#1]}}
\def\litem{\par\hang\ltextindent}

\def\ltextindent#1{\noindent\rlap{#1}\indent\ignorespaces}

\def\fine{\unskip\kern 6pt\penalty500\raise-1pt\hbox
{\vrule\vbox  to8pt{\hrule width 6pt\vfill\hrule}\vrule}\par}

\def\PP{\mathbb{P\hskip1pt}}

\def\0{\textbf{0}}

\def\OO{{\mathscr{O}}}

\def\Pic{\mathop{\rm Pic}}

\def\L{{\mathscr L}}

\def\II{{\mathscr I}}


\def\OO{\mathscr{O}}
\def\CC{\mathcal{C}}

\def\JJ{\mathcal{J}}
\def\II{\mathcal{I}}
\def\L{\mathcal{L}}
\def\P{\mathcal{P}}
\def\U{\mathcal{U}}
\def\A{\mathbb{A}}
\def\C{\mathbb{C}}

\def\PP{\mathbb{P}}
\def\Z{\mathbb{Z}}
\def\WW{\mathcal{W}}
\def\T{\mathcal{T}}
\def\Pic{\mbox{Pic}}
\author{Salvatore Giuffrida - Renato Maggioni - Riccardo Re}
\title{Multiplications of Maximal Rank in the Cohomology of $\PP^1\times\PP^1$}
\date{}
\begin{document}
\maketitle
\begin{abstract}
Let $Q=\PP^1\times\PP^1$ and let $C\subseteq Q$ be a curve of type $(a,b)$ of equation $F=0$. 
The main purpose of this paper is to analize the multiplicative structure of the bi-graded module $H^1_*\OO_Q$, in particular to prove that for any $r,s\geq 0$ the multiplication map $H^1\OO_Q(r,-s)\stackrel{F}{\longrightarrow} \OO_Q(r+a,-s+b)$ induced by $F$ 
has maximal rank for the general $C$ of type $(a,b)$. Interpretations of this problem in the contexts of multilinear algebra and differential algebra are emphasized.
\end{abstract}

\section{Introduction.}
Tensors enter in many fields of pure and applied mathematics and
the most common and useful linear maps on spaces of tensors arise
from multiplications and contractions. In this paper we address a
quite natural question about multiplication and contraction maps
of symmetric tensors. We work over the field of the complex numbers
$\mathbb{C}$ and set $V=\mathbb{C}^{m+1}$,
$W=\mathbb{C}^{n+1}$. Denote symmetric powers with $S^i$.
\begin{prob}\label{prob:1} Let
$\sigma\in S^a V\otimes S^bW$ be some fixed element and $r,t$ be
some integers with $r\geq 0$ and $t\geq b$. Consider the linear
map
\begin{equation}\label{eq:contract}S^rV\otimes
S^tW^*\stackrel{\sigma}\to S^{r+a}V\otimes
S^{t-b}W^*
\end{equation} 
defined by multiplication on the first
$r$ tensor components and contraction on the last $t$ tensor
components. Is this map of maximal rank, for $\sigma$
sufficiently general in $S^aV\otimes S^bW$?\end{prob} 
Notice that
the question is very easy if $\dim W$ or $\dim V=1$, indeed in
these cases the map (\ref{eq:contract}) is given by multiplication
of symmetric tensors or by its dual, the contraction, and it is
either injective or surjective. The first non-trivial case of
Problem \ref{prob:1} appears when $\dim V=\dim W=2$. The object of this paper
is to solve Problem \ref{prob:1} in this case. We like also to point out
 two other equivalent
formulations of Problem \ref{prob:1}. Consider the variables
$\underline{x}=(x_0,\ldots,x_m)$, $\underline{y}=(y_0,\ldots,y_n)$
and the derivations
$\underline{\partial}=(\partial{y_0},\ldots,\partial{y_n})$. 
We denote with $\C[\underline{x}]_i$ the vector space of
homogeneous polynomials of degree $i$, for any $i\geq 0$. 
\begin{prob}\label{prob:2} Consider a differential operator
$D\in\C[\underline{x}]_a\otimes\C[\underline{\partial}]_b$ and the
linear map
\begin{equation}\label{eq:deriv}
D: \C[\underline{x}]_r\otimes\C[\underline{y}]_t\to
\C[\underline{x}]_{r+a}\otimes\C[\underline{y}]_{t-b}.
\end{equation} Is this map of maximal rank if $D$ is sufficiently
general in $\C[\underline{x}]_a\otimes\C[\underline{\partial}]_b$?
\end{prob}
Now let $\PP^m$ and $\PP^n$ be projective spaces over $\C$ of
dimensions $m,\,n$, respectively, $Q=\PP^m\times\PP^n$ their
product and $\pi_1,\,\pi_2$ the first and second projection, respectively.
Recall that $\Pic(Q)\cong \Z\times\Z$, with basis
$\OO_Q(1,0)=\pi_1^*\OO_{\PP^m}(1)$ and $\OO_Q(0,1)=\pi_2^*\OO_{\PP^n}(1)$.
In these notations, one may consider the following third version
of Problem \ref{prob:1}.
\begin{prob}\label{prob:3} Consider the multiplication map
\begin{equation}\label{eq:multip}
H^n\OO_Q(r,-t-n-1)\stackrel{\sigma}\to
H^n\OO_Q(r+a,-t+b-n-1),\quad t\ge b
\end{equation} 
with $\sigma\in
H^0\OO_Q(a,b)$ a form of bi-degree $(a,b)$. Is this map of
maximal rank if $\sigma$ is sufficiently general in
$H^0\OO_Q(a,b)$? \end{prob} The fact that the three problems above
are equivalent is well known. For instance the
equivalence of Problem \ref{prob:1} and Problem \ref{prob:3} is due to the fact that
$H^n\OO_Q(r,-t-n-1)=H^0\OO_{\PP^m}(r)\otimes
H^n\OO_{\PP^n}(-t-n-1)\cong S^r V\otimes S^t W^* $, for
$V=H^0\OO_{\PP^m}(1)$ and $W=H^0\OO_{\PP^n}(1)$, by K\"unneth
formula and Serre duality. Moreover the multiplication
$H^n\OO_{\PP^n}(-t-n-1)\stackrel{\tau}\to
H^n\OO_{\PP^n}(-t+b-n-1)$, with $\tau\in H^0\OO_{\PP^n}(b)$, is
dual to the multiplication $H^0\OO_{\PP^n}(t-b)\stackrel{\tau}\to
H^0\OO_{\PP^n}(t)$. Denoting by $(y^*)^I$ the dual basis of
$y^I=y_0^{i_0}\cdots y_n^{i_n}$, with $i_0+\cdots+i_n=|I|=t$, this
map can be described as follows. For any multi-index $J$ with
$|J|=b,$ one has $(y^*)^I\cdot y^J=0$ if $J\not\subset I$ and
$(y^*)^I\cdot y^J=(y^*)^{I\setminus J}$ if $J\subset I$. This
coincides with the differentiation $D(y^*)^I$, where
$D=\partial_{y^*_0}^{i_0}\cdots\partial_{y^*_n}^{i_n} $, up to a
non zero rational number factor, showing the equivalence of
Problems \ref{prob:1} and \ref{prob:3} with Problem \ref{prob:2}. In this paper we answer
affirmatively to Problem \ref{prob:3} for $m=n=1$, using some deep facts
about the geometry of curves on the surface $Q=\PP^1\times\PP^1$.
We know also how to solve the three problems in case $a=b$ and
$n,m$ general, by a very different technique involving the
differential operator formulation of Problem \ref{prob:2}. The full solution
of Problems \ref{prob:1},\ref{prob:2},\ref{prob:3} will be the object of future investigations.
Notice that the fact that the map (\ref{eq:multip}) has maximal rank in the case of $n=m=1$ helps understanding the multiplicative stucture of the bigraded module $H^1_\ast\OO_Q$, and hence that of the Rao module of curves $C\subset Q$ when $Q$ is embedded in $\PP^3$, cfr. \cite{gm}.  
\vskip1mm
\noindent
{\bf Aknowledgements.}
The present paper is a development of the previous paper \cite{gmr} of the same 
authors. The first two authors acknowledge that its completion is mainly due to the contributions of the third author.
\section{Problem 3 for $m=n=1$. Notations and first reductions.}
Let $F\in H^0\OO_Q(a,b)$, with $a,b\geq 1$ be a form of bi-degree
$(a,b)$ and let $C$ denote the effective divisor associated to $F$;
we call it a curve of type $(a,b)$. We want to show that if
$C$ is sufficiently general in the complete linear system
$|\OO_Q(a,b)|$, then the multiplication map induced by $F$
\begin{equation}\label{eq:muF}
H^1\OO_Q(r,-t-2)\stackrel{F}\to H^1\OO_Q(r+a,-t+b-2)\end{equation}
has maximal rank. Consider the exact sequence of sheaves $$0\to
\OO_Q(r,-t-2)\stackrel{F}\to \OO_Q(r+a,-t+b-2)\to \OO_C(r+a,-t+b-2)\to 0$$
naturally associated to $F$. Since $r\geq 0$ and $t\geq b$, one
easily sees that $\ker(F)\cong H^0\OO_C(r+a,-t+b-2)$ and
$\mbox{coker}(F)\cong H^1\OO_C(r+a,-t+b-2).$ So 
Problem 3 for $m=n=1$ will be  equivalent to the following
theorem, whose proof will occupy the remainder of this paper.
\begin{thm}\label{thm:nonspecial} If $(a,b)\geq (1,1)$ and $C\in
|\OO_Q(a,b)|$ is a general curve, then for any $h\geq a$ and
$k\leq -2$ one has
\begin{equation}\label{eq:mainthm}h^0\OO_C(h,k)\cdot h^1\OO_C(h,k)=0.\end{equation}
\end{thm}
\noindent
\begin{rmk} {\rm Notice that the formulation of Problem \ref{prob:3} for $m=n=1$ in the form of Theorem \ref{thm:nonspecial} above nicely illustrates its non-triviality. For example one can easily see that any curve of type $(a,b)$ which contains a line of type $(0,1)$ fails to verify (\ref{eq:mainthm}). Even assuming smoothness, it is not difficult to find curves not verifying the statement of Theorem
\ref{thm:nonspecial} for infinitely many $(h,k)$. Consider for instance a smooth curve $E$ of type $(2,2)$ in $Q$, hence with $g(E)=1$, with $\OO_E(n,0)\cong\OO_E(0,n)$ for some $n\geq 2$. It is not difficult to show that such curves exist. Then one sees that $h^0\OO_E(ln,-ln)=\OO_E$ for any $l\geq 0$, so (\ref{eq:mainthm}) is false for $(h,k)=(ln,-ln)$}.\end{rmk}

The statement of Theorem \ref{thm:nonspecial} was proposed in our paper \cite{gmr} as Conjecture 8.1. 
Notice that the cases $a=1$ or $b=1$ are
trivial, since in these cases a general $C$ of type $(a,b)$ is
isomorphic to $\PP^1$ and the statement above is trivially
verified, since for any line bundle $L$ on $\PP^1$ one has
$h^0(L)h^1(L)=0$. So from now on we assume $(a,b)\geq (2,2)$.
Set $d=\deg\OO_C(h,k)$ and $g=(a-1)(b-1)=p_a(C)$. 
What we have to
show is that if $d+1-g\leq 0$ then $h^0\OO_C(h,k)=0$, or if
$d+1-g\geq 0$ then $h^1\OO_C(h,k)=0.$ These two problems are
equivalent. Indeed one has $\omega_C=\OO_C(a-2,b-2)$ for any curve
of type $(a,b)$ on $Q$, hence
$\omega_C\otimes\OO_C(-h,-k)=\OO_C(a-2-h,b-2-k)$. Setting
$h'=a-2-h$ and $k'=b-2-k$ by the hypothesis $h\ge a,\, k\le-2$ 
we find $h'\leq -2$ and $k'\geq b$,
moreover $h^0\OO_C(h,k)\cdot h^1\OO_C(h,k)=h^0\OO_C(h',k')\cdot
h^1\OO_C(h',k')$ hence the problem for $h,k$ is equivalent to the
one for $h',k'$, up to the automorphism of $Q=\PP^1\times\PP^1$
that interchanges the two rulings. An important final observation
is that one only needs to construct a single curve satisfying the
statement of Theorem \ref{thm:nonspecial}, since, by
semi-continuity, the statement will then hold also on a Zariski
open subset of $|\OO_Q(a,b)|$.
\subsection{Further numerical reductions for  $(h,k).$}
The proof of Theorem
\ref{thm:nonspecial} has been reduced to show that if $d+1-g\leq
0$ then $h^0\OO_C(h,k)=0$. Notice that if $(h,k)\leq (\bar{h},\bar{k})$ 
then there exists some sheaf embedding
$\OO_C(h,k)\subset \OO_C(\bar{h},\bar{k})$. So, if
$\bar{d}=\deg\OO_C(\bar{h},\bar{k})$ satisfies $\bar{d}+1-g\leq
0$, and one is able to show that $h^0\OO_C(\bar{h},\bar{k})=0$,
one immediately deduces also $h^0\OO_C(h,k)=0$. So we assume that
this reduction is no more possible, that is $d\leq g-1$,
$\deg\OO_C(h+1,k)>g-1$ and $\deg\OO_C(h,k+1)>g-1.$ This implies
\[\begin{array}{c} g-1-b<d\leq g-1\\
g-1-a<d\leq g-1,\end{array}\] that is $g-1\geq d>
ab-a-b-\min(a,b)$. Now let us write
\[\begin{array}{lr}
h=\alpha+ma &  \mbox{with }-1\leq \alpha\leq a-2\\
 k=\beta-nb & \mbox{with }
-1\leq \beta\leq b-2\end{array}\] 
Notice that $m,n>0$, due to the
limitations for $h,k$ . Assume that
$(\alpha,\beta)\not=(-1,-1)$. We want to show that in this case
one has $m=n$. We have $d=\alpha b+\beta a+(m-n)ab\leq
g-1=ab-a-b$. Since $(\alpha,\beta)>(-1,-1)$, one sees easily that
$m\leq n$. If it were $m<n$, then $d\leq \alpha b+\beta a -ab\leq
(a-2)b+(b-2)a-ab=ab-2a-2b< ab-a-b-\min(a,b)$, which is against our
assumptions on $d$. So we have $m=n$. We have proved the following
proposition.
\begin{prop}\label{prop:alphabeta} Under the assumptions above, if
$(\alpha,\beta)>(-1,-1)$ we have
$\OO_C(h,k)=\OO_C(\alpha+ma,\beta-mb)$.
 \end{prop}
 If $\alpha=\beta=-1$ we have
$\OO_C(h,k)=\OO_C(-1+ma,-1-nb)$ and $d=-b+mab-a-nab\leq
g-1=ab-a-b$, so $m-n\leq 1$. Moreover $d>ab-a-b-\min(a,b)$, so
$(m-n)ab>ab-\min(a,b)$, hence $m-n>0$. So in the case
$\alpha=\beta=-1$ we have $m=n+1$ and $\OO_C(h,k)=\OO_C(-1+ma,
b-1-mb)$. We have proved the following.
\begin{prop}\label{prop:cases-1}
Under the assumptions above, if $(\alpha,\beta)=(-1,-1)$ we have
$\OO_C(h,k)=\OO_C(-1+ma,b-1-mb)$.
 \end{prop}
In the next section we will prove Theorem \ref{thm:nonspecial} in
this last case by considering a very particular class of curves of
type $(a,b)$.
\subsection{A special class of $(a,b)$ curves}
 Take any $a$ pairwise distinct $(1,0)$ lines with equations
$L_1=0,\ldots,L_a=0$ and $b$ pairwise distinct $(0,1)$ lines with
equations $M_1=0,\ldots,M_b=0$. We denote $G$ the grid formed by the points
$P_{ij}$ determined by the equations $L_i= M_j=0$, for $i=1,\ldots,a$ and
$j=1,\ldots,b$. Then the following result holds.
\begin{lm}\label{lemma:C0} The general $(a,b)$ curve $C_0$ containing $G$ is smooth and has the property $$\OO_{C_0}(a,-b)\cong\OO_{C_0}.$$ \end{lm}
\begin{proof} Consider the $(a,0)$ curve $L=L_1\cdots L_a$ and the $(0,b)$ curve $M=M_1\cdots M_b$. We can cover $Q$ with charts $(u,v)$  parametrizing affine open sets $U\cong\A^1\times\A^1\subset Q$, in such a way that $L=(u-\lambda_1)\cdots(u-\lambda_a)=l(u)$ and $M=(v-\mu_1)\cdots(v-\mu_b)=m(v)$ on $U$, with all $\lambda$'s and $\mu$'s not equal to $0$. Then we consider a curve $C$ of type $(a,b)$, containing $G$, with equation in $U$ of the form 
\begin{equation}\label{eq:C0special} l(u)v^b-h(u)m(v)=0,\end{equation}
with $\deg h(u)=a$. Then a singular point of $C$ in $U$ must satisfy $l(u)v^b-m(v)h(u)=l'(u)v^b-m(v)h'(u)=bl(u)v^{b-1}-m'(v)h(u)=0.$ One can exclude solutions with $v=0$ by choosing $h(u)$ without multiple roots. Similarly, solutions with $m(v)=0$ are impossible since $l(u)$ has no multiple roots. Choosing $h(u)$ such that $l(u)$ and $h(u)$ have no common roots, we see that any solution $(u_0,v_0)$ of the system above must be such that the two vectors $(v_0^b, bv_0^{b-1})$ and $(m(v_0),m'(v_0))$ must be linearly dependent, hence $v_0^bm'(v_0)-bv_0^{b-1}m(v_0)=0$. So $v_0$ varies in a specified finite set $F$. Moreover $u_0$ must be a common root of $l(u)v_0^b-m(v_0)h(u)=0$ and $l'(u)v_0^b-m(v_0)h'(u)=0$, that is a multiple root of $l(u)v_0^b-m(v_0)h(u)=0$. One can exclude this possibility by choosing $h(u)$ so that $l(u)c^b-m(c)h(u)=0$ has no multiple roots for any $c\in F$. So there exists $C$ of the form (\ref{eq:C0special})  smooth on $U$. Hence the general curve of type $(a,b)$ containing $G$ is smooth on $U$, and since we can cover $Q$ with four such open affines $U$, we see such a general $C$ is smooth everywhere. \\
Finally we see that if $C_0$ is smooth and contains $G$, one has
$\OO_{C_0}(a,-b)\cong\OO_{C_0}(C_0.L-C_0.M)=\OO_{C_0}(G-G)=\OO_{C_0}.$\end{proof}
\begin{cor}\label{cor:cases-1}
The conclusion of Theorem \ref{thm:nonspecial} holds for
$$\OO_C(h,k)=\OO_C(-1+ma,b-1-mb).$$
\end{cor} 
\begin{proof} Let $C_0$ be an
$(a,b)$-curve as in the lemma above. Then
$\OO_{C_0}(a,-b)=\OO_{C_0}$, whence
$\OO_{C_0}(h,k)=\OO_{C_0}(-1,b-1)$. Then one has
$H^0\OO_{C_0}(-1,b-1)=\ker(H^1\OO_Q(-1-a,-1)\stackrel{F}\to
H^1\OO_Q(-1,b-1))=(0)$. \end{proof} 

%\begin{rmk} It is
%also true that any smooth $(a,b)$ curve $C\subset Q$
% such that $\OO_C(a,-b)=\OO_C$ contains $\infty^1$ grids like $G$ and
% has equation of the form $LH+MK=0$ with $L,K$ of type $(a,0)$ and $H,M$
% of type $(0,b)$. We omit the proof since we do not need this result in this
% paper.
%\end{rmk}

Now we are left with $\OO_C(h,k)$ as in Proposition
\ref{prop:alphabeta}. This will be the object of the remaining
sections.

\section{Completion of the proof of Theorem \ref{thm:nonspecial}.}
Given $\alpha$ and $\beta$ as in Proposition \ref{prop:alphabeta},
we set $\hat{\alpha}=a-2-\alpha$ and $\hat{\beta}=b-2-\beta$.
Notice that
$\OO_C(\hat{\alpha},\hat{\beta})=\omega_C(-\alpha,-\beta)$.\\
Let $G$ be the grid of $ab$ points in $Q$ introduced in the
preceding section. We will need the following technical result on
the existence of a subset of $G$ with particularly good properties
for our purposes.

\begin{lm}\label{lemma:Z} There exists a subset $Z\subset G$ such that
$\deg Z=(\hat{\alpha}+1)(\hat{\beta}+1)$ and
$H^0\mathscr{I}_Z(\alpha,\beta)=H^0\mathscr{I}_Z(\hat{\alpha},\hat{\beta})=0$.
\end{lm}
To prove this fact, we need the following combinatorial result.
\begin{lm}\label{lemma:bipartite} For any fixed positive integers $r,l,N$ with $N\leq rl$, there exists a bipartite graph $g$ with $r$ right vertices and $l$ left vertices, such that every right vertex has degree $\geq \lfloor N/r\rfloor$ and every left vertex has degree $\geq\lfloor N/l\rfloor$.
\end{lm}
\begin{proof} The statement is trivial if $N=1$ or if $N=rl$. One proceeds by induction on $N$. Denoting $(v_1,\ldots,v_r)$ and $(w_1,\ldots,w_l)$ the distribution of degrees at right and at left, respectively, with $\sum v_i=\sum w_j=N$, the statement consists in producing a graph with distributions of the form
\begin{eqnarray*} \underline{v}&=(v_1,\ldots,v_r)=(v+1,\ldots,v+1,v\ldots,v)\\ \underline{w}&=(w_1,\ldots,w_l)=(w+1,\ldots,w+1,w\ldots,w).
\end{eqnarray*} 
Indeed in this case one has necessarily $v=\lfloor N/r\rfloor$ and $w=\lfloor N/l\rfloor$. Now the easy proof is left to the reader.
 \end{proof}
\begin{proof}[Proof of Lemma \ref{lemma:Z}]
We set $\gamma=\max(\alpha,\hat{\alpha})$ and
$\delta=\max(\beta,\hat{\beta})$. Consider the sub-grid
$G_{\gamma,\delta}=\{P_{i,j}\ :\ 1\leq i \leq \gamma+1,\ 1\leq
j\leq \delta+1\}$. We will construct the required set $Z$ as a
subset of $G_{\gamma,\delta}$. Since we know that
$(\alpha+1)(\beta+1)\leq(\hat{\alpha}+1)(\hat{\beta}+1),$ there
are three possibilities: $(\gamma,\delta)$ equal to
$(\hat{\alpha},\hat{\beta})$ or to $(\hat{\alpha},{\beta})$ or to
$({\alpha},\hat{\beta})$. In the first case we take
$Z=G_{\gamma,\delta}$, that is the complete intersection of
$\hat{\alpha}+1$ lines of type $(1,0)$ with $\hat{\beta}+1$ lines
of type $(0,1)$. It is then clear that
$H^0\II_Z(\hat{\alpha},\hat{\beta})=H^0\II_Z(\alpha,\beta)=0$. In
the second case we set $N=(\hat{\alpha}+1)(\hat{\beta}+1)$
$r=\hat{\alpha}+1$, $l=\beta+1$ and construct a graph $g$ as in
Lemma \ref{lemma:bipartite}. Then we define $Z=\{P_{ij}\ :\
\{ij\}\in\mbox{Edges(g)}\}$. Then on any of the $\hat{\alpha}+1$
lines $L_1,\ldots,L_{\hat{\alpha}+1}$ there are $\hat{\beta}+1$
points of $Z$ and on any of the $\beta+1$ lines
$M_1,\ldots,M_{\beta+1}$ there are at least $N/(\beta+1)\geq
\alpha+1$ points of $Z$. From this it is easy to see that
$H^0\II_Z(\hat{\alpha},\hat{\beta})=H^0\II_Z(\alpha,\beta)=0$. The
third case is dealt similarly.  \end{proof}
\subsection{Avoiding the Brill-Noether locus} 
As in the preceding sections, we denote
$d=\deg\OO_C(\alpha,\beta)=\alpha b+\beta a$, $g=(a-1)(b-1)$ and assume $d\leq g-1$. Given
$G$ and $Z$ as in Lemma \ref{lemma:C0} and Lemma \ref{lemma:Z}, we
consider the linear system $S$ parametrizing the curves $C$ of type $(a,b)$
such that $G\setminus Z \subset C$. For any $\lambda\in S$ we
denote by $C_\lambda$ the corresponding curve. We also consider
the incidence variety $\CC=\{(x,\lambda)\in Q\times S \ :\ x\in
C_\lambda\}$, which defines a flat family of curves over $S$ by
means of the second projection $p:\CC\to S$. By Lemma
\ref{lemma:C0}, we know that there exists a smooth curve $C_0$
containing $G$ such that $\OO_{C_0}(a,-b)=\OO_{C_0}$. Then on a
suitable open affine neighborhood $0\in B\subset S$, the pull-back
$\CC_B\to B$ is a flat family of smooth deformations of
$C_0=p^{-1}(0)$. We want to show that for a general $\lambda\in B$
and $C_\lambda=p^{-1}(\lambda)$ one has
\[H^0\OO_{C_\lambda}(\alpha+ma,\beta-mb)=0 \quad m>0.\] For a smooth projective curve $C$ one denotes $W_d(C)$
the Brill-Noether locus
$$W_d(C)=\{L\in \Pic^d(C)\ :\ h^0(L)\not=0\}$$
We want to prove that for a general $\lambda\in B$
$$\OO_{C_\lambda}(\alpha+ma,\beta-mb)\not\in W_d(C_\lambda).$$
By the general theory of the relative Picard scheme, see for example \cite{kl}, one can associate to the family $\CC_B$ a family
$$q:\P_d\to B,$$ together with a universal line bundle $\mathcal{U}_d$ on $\P_d\times_B\CC_B$, representing the functor on the category of algebraic schemes over $B$ which associates to any $f:X\to B$ the set of equivalence classes of line bundles on $\CC_X=X\times_B\CC_B$ such that, for any closed point $x\in X$, the restriction $\L\otimes_{\OO_B}\C(x)$ is a line bundle of degree $d$ on $C_{f(x)}$, under the equivalence relation $\L\sim \L\otimes_{\OO_B}f^\ast \mathcal{N}$ for any $\mathcal{N}\in\Pic(X)$. One has
$$q^{-1}(\lambda)\cong \Pic^d(C_\lambda)$$ for any closed point $\lambda\in B$. 
We denote $O\in \P_d$ the point corresponding to the line bundle $\OO_{C_0}(\alpha,\beta)$ on $C_0$ and with $l:B\to\P_d$ the section of $q$ defined by
\begin{equation}\label{eq:l} l(\lambda)= \OO_{C_\lambda}(\alpha+ma,\beta-mb)\in \Pic^d(C_\lambda).\end{equation} The section $l$ has the property that
$l(0)=O\in \P_d$. 
 Now we consider the $d$-th symmetric power $$\CC_B^{(d)}=\CC_B\times_B\cdots\times_B\CC_B/\Sigma_d$$ relative to $B$, with $\Sigma_d$ the permutations groud on $d$ letters, and the canonical map
 $$u:\CC_B^{(d)}\to \P_d$$ which restricts fiberwise to the Abel-Jacobi maps $u_\lambda:C^{(d)}_\lambda\to\Pic^d(C_\lambda)$. It is easy to see that $\CC_B^{(d)}$ is smooth and that
 $$u^{-1}(O)=u_0^{-1}(\OO_{C_0}(\alpha,\beta))\cong |\OO_{C_0}(\alpha,\beta)|\subset C_0^{(d)}.$$
 Then we define the subvariety
$$\WW=u(\CC_B^{(d)})\subset \P_d.$$ As a set, one has $\WW=\bigcup_{\lambda}W_d(C_\lambda)$. We want to show that
$l(\lambda)\not\in \WW$ for some $\lambda\in B$. We denote by
$\T_{O}(\WW)$ the tangent cone to $\WW$ in the tangent space
$T_{O}\P_d$. Our idea is to show that the image of the differential
$dl$ is not contained in $\T_O(\WW)$, i.e. there exists some $v\in
T_0B$ such that $dl(v)\not\in\T_O(\WW)$. Indeed, if this is the
case, then for a deformation $C_\lambda$ of $C_0$ in the direction
$v$ one gets $l(\lambda)\not\in\WW$ for
general $\lambda$.
\subsection{Tangent cone to $\WW$.}
First of all, observe that by repeating the same construction of $\P_d$ in the case $d=0$ one obtains the relative Picard group $\P_0\to B$, with fibers $\Pic^0(C_\lambda)$. One obtains an isomorphism $\P_d\cong\P_0$, as schemes over $B$, by changing the universal line bundle $\U_d$ to $\U_d\otimes q^*\OO_{\CC_B}(\alpha,\beta)^{-1}$. In other words, we can uniformly identify each $\Pic^d(C_\lambda)$ with the jacobian variety $\Pic^0(C_\lambda)\cong\JJ(C_\lambda)$ by sending $\OO_{C_\lambda}(\alpha,\beta)$ to $0\in\JJ(C_\lambda)$. We keep calling $\WW$ the subvariety of $\P_0$ produced by means of the above identification, and $O\in \P_0$ the point corresponding to $\OO_{C_0}(\alpha,\beta)\in \P_d$. Notice that $O\in\P_0$ corresponds to the origin $0\in \JJ(C_0)\subset \P_0$.\\
 We recall that by Serre duality and adjunction, one can identify the tangent spaces
$T_0\JJ(C_\lambda)=H^1(\Omega^1_{C_\lambda})$ with $H^0\OO_Q(a-2,b-2)^\vee$. So the
tangent space $T_O\P_0$ is given by:
$$T_O\P_0=H^0\OO_Q(a-2,b-2)^\vee\times T_0B.$$
The following lemma describes the tangent cone to $\WW$ in $T_O\P_0$
\begin{lm}\label{lm:coneW}
Let $\T$ be the tangent cone to $W_d(C_0)$ at the point $\OO_{C_0}(\alpha,\beta)$. Then the  following facts hold.\\
1) The tangent cone $\T_O\WW$ is supported on an irreducible closed set.\\
2) As a set, $\T_O\WW=\T\times T_0B$.
\end{lm}
\begin{proof} 1) It is easy to show that the map $u:\CC_B^{(d)}\to \WW$ is birational and it satisfies all the hypotheses of Lemma 1.1 chapter II of \cite{acgh} and the subsequent corollary. The conclusion is that, as a set, $\T_O\WW$ is given as the image  of the normal bundle $\mathcal{N}$ to $u^{-1}(O)$ by means of the differential $du$. The fiber $u^{-1}(O)$ is the projective space $|\OO_{C_0}(\alpha,\beta)|$, as observed above. Hence $\T_O\WW$ is irreducible.\\
2) For any $\lambda\in B$, 
Kempf's theorem (see \cite{acgh} p. 241) describes the tangent cone at 
$\OO_{C_\lambda}(\alpha,\beta)$ of
$W_d(C_\lambda)\subset \Pic^d(C_\lambda)$ as the affine cone
$$ \bigcup_{\sigma\in H^0\OO_{C_\lambda}(\alpha,\beta)} \mu(\sigma\otimes
H^0\OO_{C_\lambda}(\hat{\alpha},\hat{\beta}))^{\perp}\subset
H^0(\OO_{C_\lambda}(a-2,b-2))^\vee,$$ where
$\mu: H^0\OO_{C_\lambda}(\alpha,\beta)\otimes
H^0\OO_{C_\lambda}(\hat{\alpha},\hat{\beta})\to
H^0\OO_{C_\lambda}(a-2,b-2)$ is the multiplication map of global
sections. Recall that one may also identify
$H^0\OO_{C_\lambda}(\alpha,\beta)\cong H^0\OO_{Q}(\alpha,\beta)$
and $H^0\OO_{C_\lambda}(\hat{\alpha},\hat{\beta})\cong
H^0\OO_{Q}(\hat{\alpha},\hat{\beta})$ and $\mu$ is identified with
the analogous multiplication map on $Q$
$$\mu_Q: H^0\OO_{Q}(\alpha,\beta)\otimes
H^0\OO_{Q}(\hat{\alpha},\hat{\beta})\to H^0\OO_{Q}(a-2,b-2).$$ As
a result, we see that the tangent cones at
$\OO_{C_\lambda}(\alpha,\beta)$ of $W_d(C_\lambda)$ are all the
same cone
$$\T=\bigcup_{\sigma\in
H^0\OO_Q(\alpha,\beta)} \mu(\sigma\otimes
H^0\OO_Q(\hat{\alpha},\hat{\beta}))^{\perp}$$
 inside the space
$H^0\OO_Q(a-2,b-2)^\vee$. By the identifications
$\Pic^d(C_\lambda)\cong J(C_\lambda)$ which map
$\OO_{C_\lambda}(\alpha,\beta)$ to $0\in J(C_\lambda)$, we see
that the tangent cones to $\WW(\lambda)$ inside
$H^0\OO_Q(a-2,b-2)\cong T_0J(C_\lambda)$ are all equal to $\T$. It
follows that the tangent cone $\T_O\WW$ contains the product cone
$\T\times T_0B$. By 1) and the
equality of dimensions, we are done.\end{proof}

\subsection{Computation of $dl$}
Let $l$ be the section of $\P_d\to B$ defined in (\ref{eq:l}).
Recall that for any curve $C$ the Picard group $\Pic^0(C)$ is
identified with the jacobian variety $J(C)$ by means of the
Abel-Jacobi map
$$D=\sum_i (P_i-Q_i)\mapsto
\sum_i\left(\int_{Q_i}^{P_i}\omega_1,\ldots,\int_{Q_i}^{P_i}\omega_g\right),$$
with $\omega_1,\ldots,\omega_g$ a basis of $H^0(\Omega^1_C)$, and
integrals defined along arbitrary paths in $C$ from $Q_i$ to
$P_i$.\\
 Now let $C_t$ be a sub-family of $\CC_B\to B$ parametrized by a line
$t\mapsto tv\in B$, for a given tangent vector $v\in T_0 B$. We know that $$l(tv)=\OO_{C_t}(ma,-mb)=\OO_{C_t}(mX_t-mY_t),$$ with $X_t=\sum_{i=1}^N P_i(t)$ and
$Y_t=\sum_{i=1}^N Q_i(t)$, $N=\deg Z$ and  $P_i(0)=Q_i(0)\in Z$ for $i=1,\ldots,N$, hence in particular $X_0=Y_0=Z$. Now let us fix our attention on one point $P=P_i(0)=Q_i(0)\in Z$ and assume that the curves $C_t$ have equations $f(u,v,t)=0$ in affine coordinates $u,v$ centered at $P$. We may assume without loss of generality that the partial derivative $f_u$ does not vanish in a neighborhood of $(u,v,t)=(0,0,0)$. Consider a basis of rational functions of $H^0\OO_Q(a-2,b-2)$ given by $g=(a-1)(b-1)$ linearly independent rational functions
$h_1(u,v),\ldots,h_g(u,v)$ such that $\mbox{div}(h_j)+(a-2)L+(b-2)M\geq 0$  for any $j=0,\ldots,g$,  with $L$ the line with equation $u=0$ and  $M$ the line with equation $v=0$. Then we know that for any $t$ in a neighborhood of $t=0$ the differential forms $(\omega_1,\ldots,\omega_g)=((h_1/f_u)dv,\ldots, (h_g/f_u)dv)$ are a basis of $H^0\Omega^1_{C_t}$. For any fixed index $i$ one can consider the Abel-Jacobi map
$$t\mapsto \left(\int_{Q_i(t)}^{P_i(t)}\omega_1,\ldots,\int_{Q_i(t)}^{P_i(t)}\omega_g\right)=s_i(t).$$
Up to suitably shrinking $B$, one can uniformly choose the integration paths $\gamma_t$ for any $t\in B$ in the following way. We can assume  $\gamma_t\subset C_t$ defined by a function $\gamma_t(s)=\gamma(s,t)=(u(s,t),v(s,t))$, with $(s,t)\in [0,1]\times B$, satisfying the following conditions:
$\gamma(s,t)$ is a $\CC^\infty$ function,
$\gamma(s,0)=P$ for any $s\in [0,1]$ and $f(\gamma(s,t),t)\equiv 0$.
Then we will use the following formula.
\begin{lm} Let $a(u,v,t)$ be a $\CC^1$ function and $\gamma(s,t)$ a family of paths satisfying the assumptions above. Then
$$\left.\frac{d}{dt}\int_{\gamma_t}adv \right| _{t=0}=a(0,0,0)(v_t(1,0)-v_t(0,0)).$$\end{lm}
\begin{proof} \[
\left.\frac{d}{dt}\int_{\gamma_t}adv \right| _{t=0}= \left.\frac{d}{dt}\int_0^1 a v_sds \right| _{t=0}
= \left.\int_0^1 ((a_u u_t+a_v v_t+a_t)v_s+av_{st})ds \right| _{t=0}.\] Now, since $\gamma(s,0)=(u(s,0),v(s,0))\equiv(0,0)$ and hence $v_s(s,0)\equiv 0$, we get
\[\left.\frac{d}{dt}\int_{\gamma_t}adv \right|_{t=0}=\int_0^1 a(\gamma(s,0),0)v_{st}(s,0))ds=a(0,0,0)(v_t(1,0)-v_t(0,0)).\]
\end{proof}
We apply this result to
$$l(tv)=\sum_{i=1}^N m\left( \int_{Q_i(t)}^{P_i(t)}\omega_1,\ldots,\int_{Q_i(t)}^{P_i(t)}\omega_g\right),$$  writing $\omega_i$ in the form $adv$ for any $i=1,\ldots,g$. We get
$$\left.\frac{dl(tv)}{dt}\right|_{t=0}= m\sum_{i=1}^N \lambda_i(\omega_1(P_i),\ldots,\omega_g(P_i))=
m\sum_{i=1}^N \lambda_i\Phi(P_i),$$
with $\Phi$ representing the canonical map $\phi:C_0\to \PP^{g-1}$ and the coefficient $\lambda_i$ equal to the $v$ component of $P_i'(t)-Q_i'(t)$ for any $i=1,\ldots,N$.
\vskip1mm\noindent
{\em Claim: One can choose the vector $v\in T_0B$ and hence the family $C_t$ in such a way that all $\lambda_i$'s above are non-zero.}
\vskip1mm\noindent
Indeed we know that the $P_i(t)$ belong to $(1,0)$ lines and the $Q_i(t)$ belong to $(0,1)$ lines. So in suitable affine coordinates $(u,v)$ one sees that $Q_i(t)=(u(t),c)$ with $c$ some constant. Hence the $v$ component of $Q_i'(t)$ is zero, so the only contribution to $\lambda_i$ comes from $P'_i(t)$. In the family of curves $\CC_B$ one can find a subfamily of curves $C_t$, with $t$ in a suitable neighborhood of $0\in \C$, passing through the points $P_i(t)=P_i+(c_i,t)$, for any $i=1,\ldots,N$, with $c_i=u(P_i)$ and $\{P_1,\ldots,P_N\}=Z$. This is possible since the dimension of the linear system of curves $S$ passing through $G\setminus Z$, of which $B$ is an open set, is $ab+a+b-(ab-N)=a+b+N>N$. Hence we find $\lambda_i=1$ for any $i=1,\ldots,N$, proving the claim.
\begin{proof}[Conclusion of the proof.]Assume we have a family $C_t$ as above, with
tangent vector $v\in T_0B$. Then $dl(v)\in T_O\P_0=H^0\OO_Q(a-2,b-2)^\ast\times T_0 B$ has the form
\[dl(v)=\left(m\sum_{i=1}^N\lambda_i\Phi(P_i),v \right).\]
By the discussion above, we see that $dl(v)\in \T_O\WW$ if and only if 
\[\sum_{i=1}^N\lambda_i\Phi(P_i)\in \T.\]
Assume that this fact holds. Then there exists $\sigma\in H^0\OO_Q(\alpha,\beta)$ such that for any $\tau\in H^0\OO_Q(\hat{\alpha},\hat{\beta})$, the hyperplane of $\PP H^0(a-2,b-2)^\vee$ represented by $\sigma\tau\in H^0\OO_Q(a-2,b-2)$ vanishes on $R=\sum_{i=1}^N\lambda_i\Phi(P_i)$. But now recall that the set $Z=\{P_1,\ldots,P_N\}$ imposes independent conditions to $H^0\OO_Q(\hat{\alpha},\hat{\beta})$, so for any $i=1,\ldots,N$ one finds a non-zero form $\tau_i\in H^0\OO_Q(\hat{\alpha},\hat{\beta})$ such that $\tau_i(P_i)\not=0$ and $\tau_i(P_j)=0$ for any $j\not=i$. So we have
 $$0=(\sigma\tau_i)(R)=\lambda_i\sigma(P_i)\tau_i(P_i),$$ from which it follows $\sigma(P_i)=0$. Hence $\sigma$ vanishes on $Z$, which contradicts the fact that $H^0\mathcal{I}_Z(\alpha,\beta)=0$.\end{proof}

\newpage
\noindent
Salvatore Giuffrida\\
{\tt giuffrida@dmi.unict.it}\\
Renato Maggioni\\
{\tt maggioni@dmi.unict.it}\\
Riccardo Re\\
{\tt riccardo@dmi.unict.it}\\
Dipartimento di Matematica e Informatica\\
Viale A. Doria, 6\\
95125 Catania, Italy.
\end{document}